\documentclass [11pt]{amsart}
\usepackage {amsmath, amssymb, amscd, graphicx, color, url, epsf}
\newtheorem {theorem}{Theorem}

\newtheorem {lemma}[theorem]{Lemma}
\newtheorem {proposition}[theorem]{Proposition}
\newtheorem {corollary}[theorem]{Corollary}

\def\zz {{\mathbb{Z}}}
\def\rr {{\mathbb{R}}}

\def\ff {{\mathbb{F}}}
\def\qq {{\mathbb{Q}}}

\def\hh {{\mathcal{H}}}
\def\tt {{\mathbb{T}}}
\def\xx {{\mathbf{x}}}
\def\yy {{\mathbf{y}}}
\def\aa {{\mathbf{a}}}
\def\uu {{\mathbf{u}}}
\def\aalpha {{\boldsymbol{\alpha}}}
\def\bbeta {{\boldsymbol{\beta}}}
\def\ggamma {{\boldsymbol{\gamma}}}

\def\del {{\partial}}
\def\per {{\Pi}}
\def\dom {{\mathcal{D}}}
\def\kh {{\widetilde{Kh}}}
\def\hfk {{\widehat{HFK}}}

\def\rk {{\operatorname{rank}}}
\def\wsigma {{\widehat {\Sigma}}}
\def\sym {{\mathrm{Sym}}}
\def\ta {{\tt_{\alpha}}}
\def\tb {{\tt_{\beta}}}
\def\tc {{\tt_{\gamma}}}
\def\td {{\tt_{\delta}}}

\def\fin {{\hfill \square}}

\def\det {{\operatorname{det}}}
\def\signature {{\operatorname{signature }}}

\def\To {{\ \rightarrow \ }}
\def\too {{\ \longrightarrow \ }}
\def\poscrossing {{L_+}}
\def\negcrossing {{L_-}}
\def\orsmoothing {{L_v}}
\def\hsmoothing {{L_h}}
\begin{document}

\title
[On the Khovanov and knot Floer homologies of quasi-alternating links]
{On the Khovanov and knot Floer homologies of quasi-alternating links}

\author [Ciprian Manolescu]{Ciprian Manolescu}
\thanks {CM was supported by a Clay Research Fellowship.}
\address {Department of Mathematics, Columbia University\\ New York, NY 10027}
\email {cm@math.columbia.edu}

\author [Peter Ozsv\'ath]{Peter Ozsv\'ath}
\thanks {PSO was supported by NSF grant numbers DMS-0505811 and FRG-0244663.}
\address {Department of Mathematics, Columbia University\\ New York, NY 10027}
\email {petero@math.columbia.edu}

\begin {abstract} 
  Quasi-alternating links are a natural generalization of alternating links. In this paper, we show that quasi-alternating links are ``homologically thin'' for
  both Khovanov homology and knot Floer homology. In particular, their
  bigraded homology groups are determined by the signature of
  the link, together with the Euler characteristic of the respective homology (i.e. the Jones or the Alexander polynomial).  The proofs use the exact
  triangles relating the homology of a link with the homologies of its
  two resolutions at a crossing.
\end {abstract}

\maketitle

\section {Introduction}

In recent years, two homological invariants for oriented links $L
\subset S^3$ have been studied extensively: Khovanov homology and knot
Floer homology. Our purpose here is to calculate these invariants for
the class of quasi-alternating links introduced in \cite{OScovers}, which generalize alternating links.

The first link invariant we will consider in this paper is {\em Khovanov's reduced homology}
(\cite{K},\cite{K2}).  This invariant takes the form of a bigraded
vector space over $\zz/2\zz$, denoted $\kh^{i,j}(L)$, whose Euler
characteristic is the Jones polynomial in the following sense:
$$ \sum_{i \in \zz, j \in \zz + \frac{l-1}{2}} (-1)^i q^j \rk \ \kh^{i,j}(L) = V_L(q),$$
where $l$ is the number of components of $L$.
The indices $i$ and $j$ are
called the homological and the Jones grading, respectively. (In our
convention $j$ is actually half the integral grading $j$ from
\cite{K}.) The indices appear as superscripts because Khovanov's
theory is conventionally defined to be a cohomology theory.  It is
also useful to consider a third grading $\delta,$ described by the
relation $\delta = j - i.$

Khovanov's original definition gives a theory whose Euler
characteristic is the Jones polynomial multiplied by the factor $q^{1/2}+q^{-1/2}$; for the reduced theory, the Euler
characteristic is the usual Jones polynomial, i.e. normalized so that it takes the value $1$ on the unknot, cf.~\cite{K2}.  Note that $\kh$
can be also be defined with integer coefficients, but then it depends
on the choice of a component of the link.  Nevertheless, $\kh$ is a link invariant over $\zz/2\zz$; see \cite[Section 5]{OScovers} or \cite[Section 3]{Sh}.

The other homological link invariant that we consider in this paper  is {\em knot Floer homology}. This theory was independently introduced by Szab\'o and the second author in \cite{OSknots},
and by Rasmussen \cite{R1}.  In its
simplest form, it is a bigraded Abelian group $\hfk_i(L, j)$ whose
Euler characteristic is (up to a factor) the Alexander-Conway
polynomial $\Delta_L(q)$:
$$ \sum_{j \in \zz, i \in \zz+\frac{l-1}{2}} 
(-1)^{i+\frac{l-1}{2}} q^j \rk \ \hfk_i(L,j) = (q^{-1/2}-
q^{1/2})^{l-1}\cdot \Delta_L(q).$$ Knot Floer homology was originally
defined using pseudo-holomorphic curves, but there are now also
several combinatorial formulations available, cf. \cite{MOS},
\cite{MOST}, \cite{SW}, \cite{OSsingular}.  The two gradings $i$ and
$j$ are called the {\em Maslov} and {\em Alexander gradings}
respectively; we also set $\delta=j-i$ as above. Knot Floer homology 
detects the genus of
a knot \cite{OSgenus}, as well as whether a knot is fibered \cite{N}.
There exists also an improvement, called link Floer homology
(\cite{OSlinks}, \cite{OSnorm}), which detects the Thurston norm of
the link complement, but that theory will not be discussed in this
paper. Also, even though $\hfk$ can be defined with integer coefficients, in this paper we will only consider it with coefficients in the field $\ff = \zz/2\zz.$

For many classes of links (including most knots with small crossing
number), the Khovanov and knot Floer homologies over $R=\zz$ or $\ff$
take a particularly simple form: they are free $R$-modules supported
in only one $\delta$-grading. We call such links {\it Khovanov
homologically thin (over $R$)}, or {\it Floer homologically thin (over
$R$)}, depending on which theory we refer to. Various versions of
these definitions appeared in \cite{BN},
\cite{R1}, \cite{K2}, \cite{R2}.  Further, it turns out that typically
the $\delta$-grading in which the homology groups are supported equals
$-\sigma/2,$ where $\sigma$ is the signature of the link. When this is
the case, we say that the link is (Khovanov or Floer) {\it
homologically $\sigma$-thin}. (Floer homologically $\sigma$-thin knots
were called {\em perfect} in \cite{R0}.)

If a link $L$ is homologically $\sigma$-thin over $R=\zz$ or $\ff$ for
a bigraded theory $\hh$ (where $\hh$ could denote either $\kh$ or
$\hfk$), then $\hh(L)$ is completely determined by the signature
$\sigma$ of $L$ and the Euler characteristic $P(q)$ of $\hh$ (the
latter being either the Jones or a multiple of the Alexander
polynomial). Indeed, if $P(q) = \sum a_j q^j,$ we must have: $$
\hh^{i,j}(L) \simeq
\begin{cases}
R^{|a_j|} & \text{ if } i = j+ \frac{\sigma}{2} \\
0 & \text{otherwise}. 
\end {cases}
$$

In the world of Khovanov homology, the fact that the vast majority
(238) of the 250 prime knots with up to 10 crossings are homologically
$\sigma$-thin was first observed by Bar-Natan, based on his
calculations in \cite{BN}. Lee \cite{L} showed that alternating links
are Khovanov homologically $\sigma$-thin. Since 197 of the prime knots
with up to 10 crossings are alternating, this provides a partial
explanation for Bar-Natan's observation.

At roughly the same time, a similar story unfolded for knot Floer
homology. Rasmussen \cite{R0} showed that 2-bridge knots are Floer
homologically $\sigma$-thin; and this result was generalized in
\cite{OSalt} to all alternating knots.

\begin {figure}
\begin {center}
\input {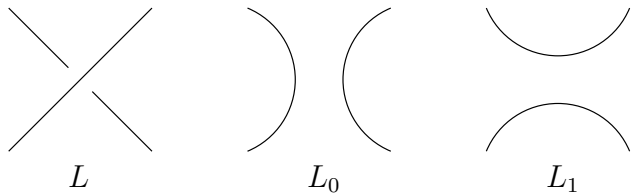}
\caption {The links in the unoriented skein relation.}
\label {fig:skein}
\end {center}
\end {figure}

In this paper we generalize these results to a larger
class of links, the {\em quasi-alternating links} of \cite{OScovers}.
Precisely, $\mathcal{Q}$ is the smallest set of links satisfying the
following properties:
\begin {itemize}
\item {The unknot is in $\mathcal{Q}$;}
\item {If $L$ is a link which admits a projection with a crossing such
that
\begin {enumerate}
\item {both resolutions $L_0$ and $L_1$ at that crossing (as in Figure~\ref{fig:skein}) are in $\mathcal{Q},$}
\item {$\det(L) = \det(L_0) + \det(L_1),$}
\end {enumerate}
then $L$ is in $\mathcal{Q}.$}
\end {itemize}

The elements of $\mathcal{Q}$ are called quasi-alternating links. It
it is easy to see (cf.~\cite[Lemma 3.2]{OScovers}) that alternating
links are quasi-alternating.

In this paper we prove the following:

\begin {theorem}
\label {khovanov}
Quasi-alternating links are Khovanov homologically $\sigma$-thin (over $\zz$).
\end {theorem}

\begin {theorem}
\label {floer}
Quasi-alternating links are Floer homologically $\sigma$-thin (over
$\zz/2\zz$).
\end {theorem}

For knots with up to nine crossings, Theorem~\ref{khovanov} and Theorem~\ref{floer} provide an almost complete explanation for the prevalence of homological $\sigma$-thinness (over the respective coefficient ring). Indeed, among the 85 prime knots with up to nine crossings, only two ($8_{19}$ and 
$9_{42}$) are not Khovanov homologically $\sigma$-thin, and these are also 
the only ones which are not Floer homologically $\sigma$-thin. By the results 
of \cite{OScovers}, \cite{M} and \cite{B}, 82 of the 83 remaining knots are 
quasi-alternating. (Among them, 74 are alternating.) This leaves only the knot $9_{46}$, which the authors do not know if it is quasi-alternating.

In general it is difficult to decide whether a larger, homologically $\sigma$-thin knot is quasi-alternating. It remains a challenge to find
homologically $\sigma$-thin knots that are not quasi-alternating; $9_{46}$ could be the first potential example.

A few words are in order about the strategy of proof and the
organization of the paper. Both Theorem~\ref{khovanov} and
Theorem~\ref{floer} are consequences of the unoriented skein exact
triangles satisfied by the respective theories. For Khovanov homology,
this exact triangle (which relates the homology of $L$ to that of its
resolutions $L_0$ and $L_1$, cf. Figure~\ref{fig:skein}), is immediate
from the definition of the homology groups. The only new ingredient
used in the proof of Theorem~\ref{khovanov} is an observation about
relating the gradings to the signature. We explain this in
Section~\ref{sec:khovanov} of the paper.  In fact, the proof of
Theorem~\ref{khovanov} is an adaptation of the proof of the
corresponding fact for alternating links due to Lee~\cite{L}.

For knot Floer homology, an unoriented skein exact triangle was
described by the first author in
\cite{M}. In that paper, the maps in the triangle were ungraded. In
Section~\ref{sec:floer}, we show that they actually respect the
$\delta$-grading, up to a well-determined shift. This will imply
Theorem~\ref{floer}. It is interesting to note that this strategy is
quite different from the earlier proofs for two-bridge and alternating
links,~\cite{R1}, \cite{OSalt}.

We remark that Theorem~\ref{floer} has a number of formal
consequences. The full version of knot Floer homology is a graded,
filtered chain complex over the polynomial algebra $\ff[U].$ It was
shown in \cite[Theorem 1.4 and the remark immediately after]{OSalt}
that for Floer homologically $\sigma$-thin knots, their full complex
(up to equivalence) is determined by their Alexander polynomial and
signature. Theorem~\ref{floer} implies then that this is true for
quasi-alternating knots. Furthermore, according to \cite{OSinteger}
and \cite{OSrational}, the full knot Floer complex has enough
information to determine the Heegaard Floer homology of any Dehn
surgery on that knot. Thus, the Floer homologies (over $\ff$) of Dehn
surgeries on quasi-alternating knots are determined by the Alexander
polynomial, the signature, and the surgery coefficient; we refer to
\cite{OSalt},
\cite{OSinteger}, \cite{OSrational} for the precise statements.

It is natural to expect Theorem~\ref{floer} to hold also over $\zz$.
Note that Theorem~\ref{floer}, combined with the universal
coefficients theorem, implies that quasi-alternating links are Floer
homologically $\sigma$-thin over $\qq$.

\medskip \noindent {\bf Acknowledgments.} We would like to thank John Baldwin, Matthew Hedden, Robert Lipshitz, Jacob Rasmussen and Sucharit Sarkar for helpful conversations. 

\section {The exact triangle for Khovanov homology}
\label {sec:khovanov}

\subsection {The Gordon-Litherland formula}
Let us review the definition of the Goeritz matrix, as well as the Gordon-Litherland formula for the signature, following \cite{GL}.

\begin {figure}
\begin {center}
\input {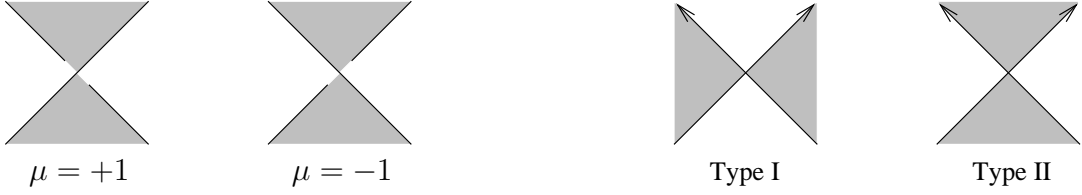}
\caption {Incidence numbers and types of crossings.}
\label {fig:gl}
\end {center}
\end {figure}

Consider an oriented link $L$ in $S^3$ with a regular, planar projection, and let $D$ be the corresponding planar diagram. The complement of the projection in $\rr^2$ has a number of connected components, which we call {\it regions}. We color them in black and white in checkerboard fashion. Let $R_0, R_1, \dots, R_n$ be the white regions. Assume that each crossing is incident to two distinct white regions. To each crossing $c$ we assign an incidence number $\mu(c),$ as well as a type (I or II), as in Figure~\ref{fig:gl}. Note that the sign of the crossing is determined by its incidence number and type.

Set
$$ \mu(D) = \sum_{c \text { of type II}} \mu(c).$$

The Goeritz matrix $G=G(D)$ of the diagram $D$ is defined as follows. For any $i, j \in \{0, 1, \dots, n\}$ with $i \neq j,$ let
$$ g_{ij} = - \sum_{c \in \bar R_i \cap \bar R_j} \mu(c).$$
Set also
$$ g_{ii} = - \sum_{i \neq j} g_{ij}.$$

Then $G$ is the $n \times n$ symmetric matrix with entries $g_{ij},$ for $i, j \in \{1, \dots, n\}.$ 

Gordon and Litherland showed that the signature of $L$ is given by the formula
\begin {equation}
\label {eq:gl}
 \sigma(L) = \signature(G) - \mu(D).
\end {equation}

(We use the convention that the signature of the right-handed trefoil is $-2.$)
Also, the determinant $\det(L)$ of a link $L$ can be defined as the non-negative integer
$$ \det(L) = | \det(G) |.$$

\subsection {The signature of resolutions.}
\begin {figure}
\begin {center}
\input {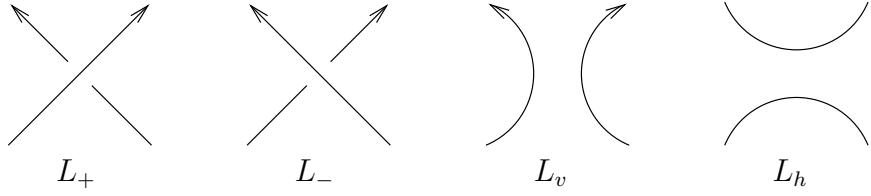}
\caption {Two possible crossings and their resolutions.}
\label {fig:fourlinks}
\end {center}
\end {figure}

Let $L \subset S^3$ be an oriented link with a fixed planar projection as before. Fix now a crossing $c_0$ in the corresponding planar diagram. If the crossing is positive (resp. negative), we set $L_+ = L$ (resp. $L_- = L$) and let $L_-$ (resp. $L_+$) be the link obtained form $L$ by changing the sign of the crossing. Further, we denote by $L_v$ and $L_h$ the oriented and unoriented resolutions of $L$ at that crossing, cf. Figure~\ref{fig:fourlinks}. (We choose an arbitrary orientation for $L_h$.) To make the connection with Figure~\ref{fig:skein}, note that if $L= L_+, $ then $L_0 = L_v$ and $L_1 = L_h,$ while if $L=L_-,$ then $L_0 = L_h$ and $L_1 = L_v.$

Denote by $D_+, D_v, D_h$ the planar diagrams of $L_+, L_v, L_h,$ respectively, differing from each other only at the chosen crossing $c_0$. 

The first equality in the lemma below (without the sign) is due to Murasugi \cite{Mur1}; the second is also inspired by a result of Murasugi from \cite{Mur2}.

\begin {lemma}
\label {signatures}
Suppose that $\det(L_v), \det(L_h) > 0$ and $\det(\poscrossing)= \det(\orsmoothing) + \det(\hsmoothing).$ Then:
$$ \sigma(\orsmoothing) - \sigma(\poscrossing) = 1$$
and
$$ \sigma(\hsmoothing)-\sigma(\poscrossing) = e,$$
where $e$ denotes the difference between the number of negative crossings in $D_h$ and the number of such crossings in $D_+.$
\end {lemma}

\begin {figure}
\begin {center}
\input {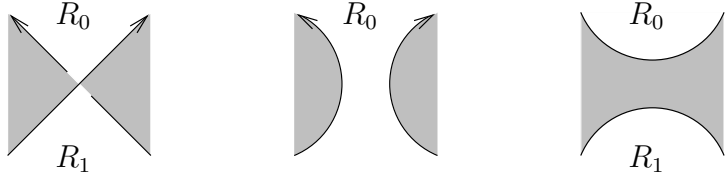}
\caption {Coloring convention at the chosen crossing.}
\label {fig:coloring}
\end {center}
\end {figure}

\noindent\textbf{Proof.} Construct the Goeritz matrices $G_+ = G(D_+), G_v = G(D_v)$ and $G_h = G(D_h)$ in such a way that $c_0$ is of Type I (and incidence number $-1$) in $D_+$, and the white region $R_0$ (the one not appearing in the Goeritz matrix) is as in Figure~\ref{fig:coloring}.

Observe now that $G_+$ and $G_h$ are bordered matrices of $G_v.$ More precisely, if $G_v$ is an $n \times n$ symmetric matrix, then there exists $a \in \rr$ and $v = (v_1, \dots, v_n)\in \rr^n$ such that
$$ G_+ = \begin {pmatrix} a & v \\ v^{T} & G_v \end {pmatrix}; \ \ \ G_h = \begin {pmatrix} a+1 & v \\ v^{T} & G_v \end {pmatrix}.$$
  
Without loss of generality (after an orthonormal change of basis), we can assume that $G_v$ is diagonal, with diagonal entries $\alpha_1, \dots, \alpha_n$. Note that these are nonzero because $\det(L_v) = |\det(G_v)| \neq 0.$

The bilinear form associated to $G_+$ can be written as
$$ aX^2 + 2\sum_{i=1}^n v_i XX_i +\sum_{i=1}^n \alpha_i X_i^ 2,$$
or
$$ \Bigl(a - \sum_{i=1}^n \frac{v_i^2}{\alpha_i}\Bigr) X^2 + \sum_{i=1}^n \alpha_i \Bigl(X_i + \frac{v_i}{\alpha_i}X\Bigr)^2.$$

A similar formula holds for the form of $G_h$, but with $a$ replaced by $a+1.$ 

If we set 
$$\beta =  a - \sum_{i=1}^n \frac{v_i^2}{\alpha_i},$$ 
then 
$$ \det(G_+) = \beta \cdot \det(G_v), \ \ \det(G_h) = (\beta + 1) \cdot \det(G_v).$$

By the condition on the determinants in the hypothesis, $|\beta| = |\beta + 1| + 1,$ so we must have $\beta < -1.$ Therefore, when we diagonalize the bilinear forms, for $G_+$ (resp. $G_h$) we get one additional negative coefficient ($\beta$, resp. $\beta + 1$) as compared to $G_v.$ Thus,
\begin {equation}
\label {signas}
 \signature(G_+) = \signature(G_h) = \signature(G_v) - 1.
\end {equation}

Since $c_0$ is of Type I, we also have $\mu(D_+) = \mu(D_v).$ Together with the Gordon-Litherland formula \eqref{eq:gl}, these identities imply
$$ \sigma(L_+) = \sigma(L_v) - 1.$$

Next, observe that when we change the direction of an arc at a crossing, both the sign and the type of the crossing are reversed, but the incidence number remains the same. If we denote by $k(\mu,t)$ the number of crossings of incidence number $\mu\in \{\pm 1\}$ and type $t \in \{I, II \}$ in $D_+$ (excluding $c_0$) which change type (and sign) in $D_h$, then
$$ \mu(D_h) - \mu(D_+) = k(+1, I) - k(-1, I) - k(+1, II) + k(-1, II).$$
This equals
$$ - \bigl( k(-1, I) + k(+1, II) \bigr) +  \bigl( k(+1, I) + k(-1, II) \bigr) = -e.$$
Using \eqref{eq:gl} and \eqref{signas} again, we get 
$$ \sigma(L_+) = \sigma(L_h) + e,$$
as desired. $\hfill \fin$

\subsection {An unoriented skein exact triangle.}
The following proposition is a simple consequence of the definition of Khovanov cohomology. It is implicit in \cite{K}, and also appeared in Viro's work \cite{Viro}. The statement below, with the precise gradings, is taken from Rasmussen's review \cite[Proposition 4.2]{R2}. It is written there in terms of Khovanov's unreduced homology, but it works just as well for the reduced version $\kh$, which we use in this paper. We work over $\zz$, so to define the reduced homology we need to mark a component for each link appearing in the triangle; we do this by marking the same point on their diagrams, away from the crossing where the links differ.

\begin{proposition}
\label{KhSkein}
(Khovanov, Viro, Rasmussen)
There are long exact sequences
$$
\cdots \to \kh^{i-e-1, j-\frac{3e}{2}-1 }(\hsmoothing) \to \kh^{i,j}(\poscrossing)
\to \kh^{i,j-\frac{1}{2}} (\orsmoothing) \to \kh^{i-e, j-\frac{3e}{2}-1 }(\hsmoothing) \to \cdots
$$
and
$$
\cdots \to \kh^{i, j+\frac{1}{2} }(\orsmoothing) \to \kh^{i,j}(\negcrossing)
\to \kh^{i-e+1,j-\frac{3e}{2}+1} (\hsmoothing) \to \kh^{i+1, j+\frac{1}{2} }(\orsmoothing) \to \cdots
$$
where $e$ is as in the statement of Lemma~\ref{signatures}.
\end{proposition}

If we forget about $i$ and $j$ and just keep the grading $\delta = j-i,$ the two triangles become
\begin {equation}
\label {1kh}
\cdots \to \kh^{*-\frac{e}{2}}(\hsmoothing) \to \kh^*(\poscrossing)
\to \kh^{*-\frac{1}{2}} (\orsmoothing) \to\kh^{*-\frac{e}{2}-1}(\hsmoothing) \to  \cdots
\end {equation}
and
\begin {equation}
\label {2kh}
\cdots \to \kh^{*+\frac{1}{2}}(\orsmoothing) \to \kh^*(\negcrossing)
\to \kh^{*-\frac{e}{2}} (\hsmoothing) \to \kh^{*-\frac{1}{2}}(\orsmoothing)\to \cdots
\end {equation}

\begin {proposition}
\label {newtriangle}
Let $L$ be a link and $L_0, L_1$ its two resolutions at a crossing as in Figure~\ref{fig:skein}. Assume that $\det(L_0), \det(L_1) > 0$ and $\det(L) = \det(L_0) + \det(L_1).$  Then there is an exact triangle:
$$ \cdots \to \kh^{*-\frac{\sigma(L_1)}{2}}(L_1) \to \kh^{*-\frac{\sigma(L)}{2}}(L) \to \kh^{*-\frac{\sigma(L_0)}{2}} (L_0) \to \kh^{*-\frac{\sigma(L_1)}{2}-1}(L_1) \to \cdots $$
\end {proposition}

\noindent\textbf{Proof.} When the given crossing in $L$ is positive, this is a re-writing of the triangle \eqref{1kh}, taking into account the result of Lemma~\ref{signatures}. Note that when following three consecutive maps in the triangle the grading decreases by one; thus, the grading change for the map between the homologies of the two resolutions is determined by the grading change for the other two maps.

The case when the crossing is negative is similar. $\hfill \fin$
\medskip

\noindent\textbf{Proof of Theorem~\ref{khovanov}.} Note that any quasi-alternating link has nonzero determinant; this follows easily from the definition. The desired result is then a consequence of Proposition~\ref{newtriangle}: the unknot is homologically $\sigma$-thin and, because of the exact triangle, if $L_0$ and $L_1$ are homologically $\sigma$-thin, then so is $L. \ \hfill \fin$

\section {The exact triangle for knot Floer homology}
\label {sec:floer}

In this section we assume that the reader is familiar with the basics
of knot Floer homology (including the version with several
basepoints), cf. \cite{OSknots}, \cite{R1}, \cite{OSlinks},
\cite{MOS}. Throughout this section we will work with coefficients in the field $\ff = \zz/2\zz$. 

\subsection {Heegaard diagrams and periodic domains}
\label {sec:periodic}

We start with a few generalities about periodic domains in Heegaard diagrams. Our discussion is very similar to the ones in \cite[Section 2.4]{HolDisk} and \cite[Section 3.4]{OSlinks}, except that here we do not ask for the periodic domains to avoid any basepoints.

Let $\Sigma$ be a Riemann surface of genus $g$. A collection $\aalpha = (\alpha_1, \dots, \alpha_n)$ of disjoint, simple closed curves on $\Sigma$ is called {\em good} if the span $S_{\alpha}$ of the classes $[\alpha_i]$ in $H_1(\Sigma; \zz)$ is $g$-dimensional. If $\aalpha$ is such a collection, we view (the closures of) the components of $\Sigma - (\cup \alpha_i)$ as two-chains on $\Sigma$ and denote by $\per_{\alpha}$ their span. Note that $\per_{\alpha}$ is a free Abelian group of rank $m=n-g+1.$ 

A {\em Heegaard diagram} $(\Sigma , \aalpha, \bbeta)$ consists of a Riemann surface $\Sigma$ together with two good collections of curves $\aalpha = (\alpha_1, \dots, \alpha_n)$ and $\bbeta = (\beta_1, \dots, \beta_n).$ (A Heegaard diagram describes a $3$-manifold $Y$; see for example \cite[Section 3.1]{OSlinks}.) We define a {\em periodic domain} in the Heegaard diagram $(\Sigma, \aalpha, \bbeta)$ to be a two-chain on $\Sigma$ that is a linear combination of the components of $\Sigma - (\cup \alpha_i) - (\cup \beta_i)$, and with the property that its boundary is a linear combination of the alpha and beta curves. (This is a slight modification of \cite[Definition 2.14]{HolDisk}.) The group of periodic domains is denoted $\per_{\alpha, \beta}.$ Let also $S_{\alpha, \beta}= S_{\alpha} + S_{\beta}$ be the span of all the alpha and beta curves in $H_1(\Sigma; \zz).$ 

\begin {lemma} 
\label {perab}
The group $\per_{\alpha, \beta}$ of periodic domains is free Abelian of rank equal to $2n+1 - \rk (S_{\alpha, \beta}).$
\end {lemma}

\noindent\textbf{Proof.} There is a map
$$ \psi_{\alpha, \beta}: \zz^{2n} \to S_{\alpha, \beta}$$
taking the first $n$ standard generators of $\zz^{2n}$ to the classes $[\alpha_i], i = 1, \dots, n,$ and the remaining $n$ standard generators to the classes $[\beta_i], i =1, \dots, n.$ There is a short exact sequence
\begin {equation}
\label {psiab}
 0 \too \zz \too \per_{\alpha, \beta} \too \ker(\psi_{\alpha, \beta}) \too 0.
\end {equation}

Indeed, the map $\per_{\alpha, \beta} \to \ker(\psi_{\alpha, \beta})$ takes a periodic domain $\dom$ to the coefficients of the alpha and beta curves appearing in $\del \dom.$ It is surjective, and its kernel is generated by the Heegaard surface $\Sigma$ itself. 

The conclusion follows immediately from the short exact sequence. $\hfill \fin$

\medskip

Note that we can view $\per_{\alpha}$ and $\per_{\beta}$ as subgroups of $\per_{\alpha, \beta}.$ Their intersection is generated by the two-chain $\Sigma.$ Therefore,
$$ \rk (\per_{\alpha} + \per_{\beta}) = 2n-1.$$

More precisely, if we denote by $S_{\alpha} \oplus S_{\beta} \cong \zz^{2g}$ the exterior direct sum, there is a short exact sequence analogous to (\ref{psiab}):
\begin {equation}
\label {psi2}
 0 \To \zz \too \per_{\alpha} + \per_{\beta} \too \ker(\zz^{2n} \to S_{\alpha} \oplus S_{\beta}) \To 0.
\end {equation}

\begin {corollary}
\label {pery}
If $S_{\alpha, \beta} = H_1(\Sigma; \zz),$ then $\per_{\alpha, \beta} = \per_{\alpha} + \per_{\beta}.$
\end {corollary}

\noindent\textbf{Proof.} The exact sequences \eqref{psiab} and \eqref{psi2} fit into a commutative diagram
$$\begin {CD}
0 @>>> \zz @>>> \per_{\alpha} + \per_{\beta} @>>> \ker(\zz^{2n} \to S_{\alpha} \oplus S_{\beta}) @>>> 0 \\
@.  @V{\cong}VV @VVV @VVV @. \\
0 @>>> \zz @>>> \per_{\alpha, \beta} @>>> \ker(\zz^{2n} \to S_{\alpha, \beta}) @>>> 0 \\
\end {CD}$$

To show that the middle vertical arrow is an isomorphism it suffices to show that the right vertical arrow is. The map $\psi_{\alpha, \beta} : \zz^{2n} \to S_{\alpha, \beta}$ factors through $ S_{\alpha} \oplus S_{\beta}$. Consider the sequence of maps
$$ \zz^{2g} \cong S_{\alpha} \oplus S_{\beta} \too S_{\alpha, \beta} \ \hookrightarrow H_1(\Sigma; \zz) \cong \zz^{2g}.$$

The hypothesis says that the last inclusion is an isomorphism, which means that the composition is surjective. Since its domain and target are both $\zz^{2g}$, the map must be an isomorphism. This shows that $S_{\alpha} \oplus S_{\beta} \too S_{\alpha, \beta}$ is an isomorphism as well.$\hfill \fin$
\medskip

Finally, a {\em triple Heegaard diagram}  $(\Sigma, \aalpha, \bbeta, \ggamma)$ consists of a Riemann surface $\Sigma$ together with three good collections of curves $ \aalpha =( \alpha_1, \dots, \alpha_n), \bbeta=(\beta_1, \dots, \beta_n), \ggamma=( \gamma_1, \dots, \gamma_n)$. A {\em triply periodic domain} is then a two-chain on $\Sigma$ that is a linear combination of the components of $\Sigma - (\cup \alpha_i) - (\cup \beta_i) - (\cup \gamma_i)$, and with the property that its boundary is a linear combination of the alpha, beta, and gamma curves. 

The group of triply periodic domains is denoted $\per_{\alpha, \beta, \gamma}.$ Set $S_{\alpha, \beta, \gamma} = S_{\alpha} + S_{\beta} + S_{\gamma} \subset H_1(\Sigma; \zz).$ A straightforward analog of Lemma~\ref{perab} then says that $\per_{\alpha, \beta, \gamma}$ is a free Abelian group of rank equal to $3n+1 - \rk (S_{\alpha, \beta, \gamma}).$

\subsection {The ungraded triangle}
\label{sec:heegaard}
The following theorem was proved in \cite{M}:

\begin {theorem}
\label {main}
Let $L$ be a link in $S^3$, and $L_0$ and  
$L_1$ the two resolutions of $L$ at a crossing, as in  
Figure~\ref{fig:skein}. Denote by $l,l_0, l_1$ the number of 
components of the links $L, L_0,$ and $L_1,$ respectively, and set $m = 
\max \{ l, l_0, l_1\}.$ Then, there is an exact triangle
$$
 \hfk(L) \otimes V^{m-l} \To \hfk(L_0)  \otimes V^{m-l_0} \To \hfk(L_1) 
 \otimes V^{m-l_1} \To \hfk(L) \otimes V^{m-l}, $$
where $V$ denotes a two-dimensional vector space over $\ff.$
\end {theorem}

Our goal will be to study how the maps in the exact triangle behave with respect to the $\delta$-grading. In order to do this, we recall how the maps were constructed in \cite{M}.

\begin {figure}
\begin {center}
\input {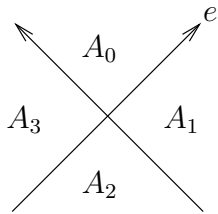}
\caption {The regions near the crossing $c_0.$ Since $c_0$ can be either negative or positive, we have not marked which strand is the overpass.}
\label {fig:czero}
\end {center}
\end {figure}

The starting point is a special Heegaard diagram which we associate to
a regular, connected, planar projection $D$ of the link $L.$ (This is
a suitable stabilization of the diagram considered in~\cite{OSalt}.)
We assume that one of the crossings in $D$ is $c_0,$ such that the two
resolutions at $c_0$ are diagrams $D_0$ and $D_1$ for $L_0$ and $L_1$,
respectively. If $D$ has $k$ crossings, then it splits the plane into
$k+2$ regions. Let $A_0, A_1, A_2, A_3$ be the regions near $c_0$ in
clockwise order, as in Figure~\ref{fig:czero}, and $e$ the edge
separating $A_0$ from $A_1.$ We can assume that $A_0$ is the unbounded
region in $\rr^2- D$. Denote the remaining regions by $A_4, \dots,
A_{k+1}.$ Let $p$ be a point on the edge $e.$ If $m =\max \{ l, l_0,
l_1\}$ is as in the statement of Theorem~\ref{main}, then we can
choose $p_1, \dots, p_{m-1}$ to be a collection of points in the plane,
distinct from the crossings and such that for every component of any
of the links $L$, $L_0$ and $L_1$, the projection of that component
contains at least one of the points $p_i$ or $p.$

\begin {figure}
\begin {center}
\begin{picture}(0,0)%
\includegraphics{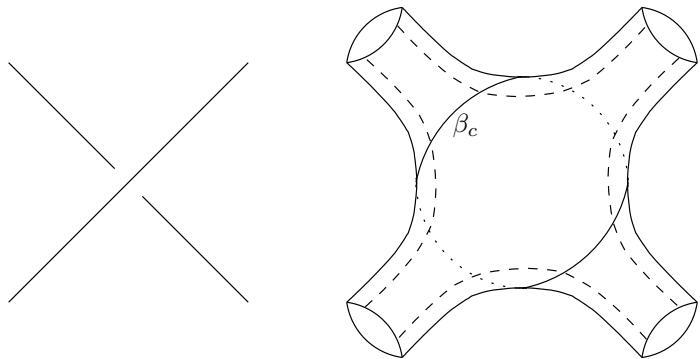}%
\end{picture}%
\setlength{\unitlength}{3947sp}%
\begingroup\makeatletter\ifx\SetFigFont\undefined%
\gdef\SetFigFont#1#2#3#4#5{%
  \reset@font\fontsize{#1}{#2pt}%
  \fontfamily{#3}\fontseries{#4}\fontshape{#5}%
  \selectfont}%
\fi\endgroup%
\begin{picture}(4349,2224)(-2261,-1973)
\put(541,-543){\makebox(0,0)[lb]{\smash{{\SetFigFont{10}{12.0}{\rmdefault}{\mddefault}{\updefault}{\color[rgb]{0,0,0}$\beta_c$}%
}}}}
\end{picture}%

\caption {Piece of the Heegaard surface $\Sigma$ associated to a crossing $c$. It contains four (or fewer) bits of alpha curves, shown in dashed lines, and  
one beta curve $\beta_c.$}
\label {fig:inters}
\end {center}
\end {figure}

\begin {figure}
\begin {center}
\input {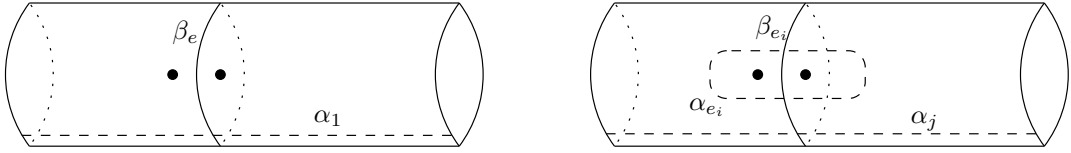}
\caption {A neighborhood of the distinguished edge $e$ (left) and a  
ladybug around some edge $e_i$ marked by $p_i$ (right).}
\label {fig:lady}
\end {center}
\end {figure}

We denote by $\Sigma$ the boundary of a regular neighborhood of $D$ in
$S^3,$ a surface of genus $g = k+1.$ To every region $A_r \ (r>0)$ we
associate a curve $\alpha_r$ on $\Sigma,$ following the boundary of
$A_r.$ To each crossing $c$ in $D$ we associate a curve $\beta_c$ on
$\Sigma$ as indicated in Figure~\ref{fig:inters}. In addition, we
introduce an extra curve $\beta_e$ which is the meridian of the knot,
supported in a neighborhood of the distinguished edge $e.$ We also
mark the surface $\Sigma$ with two basepoints, one on each side of
$\beta_e$, as shown on the left side of Figure~\ref{fig:lady}.

Furthermore, for every edge $e_i$ of $D$ containing one of the points $p_i, \ i=1, \dots, m-1,$ we introduce a ladybug,
i.e. an additional pair of alpha-beta curves on $\Sigma$, as well as an
additional pair of basepoints. This type of configuration is shown on the
right side of Figure~\ref{fig:lady}.

The surface $\Sigma$, together with the collections of alpha curves,
beta curves and basepoints, forms a multi-pointed Heegaard diagram for
$S^3$ compatible with $L$, in the sense of \cite[Definition 2.1]{MOS}.
We denote the alpha and the beta curves in the diagram by $\alpha_i,
\beta_i$ with $i = 1, \dots, n,$ where $n = g+m-1.$ We reserve the index
$n$ for the beta curve $\beta=\beta_n$ associated to the crossing
$c_0.$ Also, we let $\wsigma$ denote the complement of the basepoints
in the surface $\Sigma.$

\begin {figure}
\begin {center}
\input {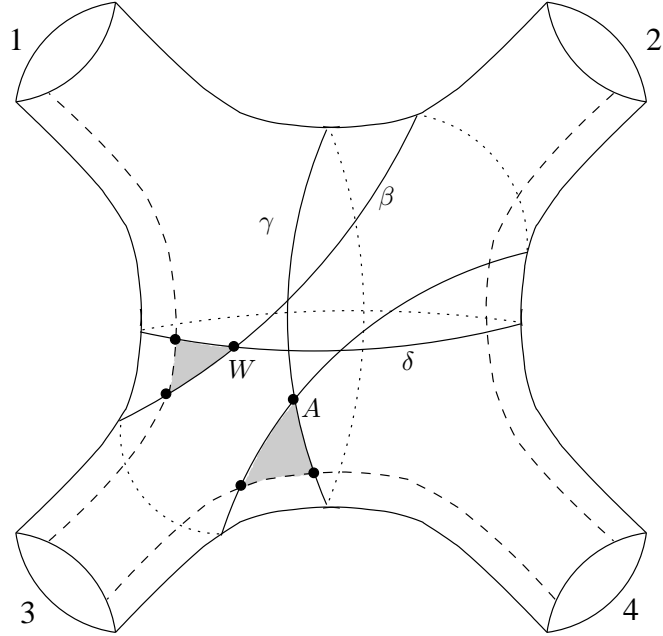}
\caption {Piece of $\Sigma$ near the crossing $c_0.$ There are three bits of alpha curves, shown dashed. This piece is joined to the rest of the diagram by four tubes, which we mark by the numbers 1,2,3,4.}
\label {fig:trio}
\end {center}
\end {figure}

\begin {figure}
\begin {center}
\input {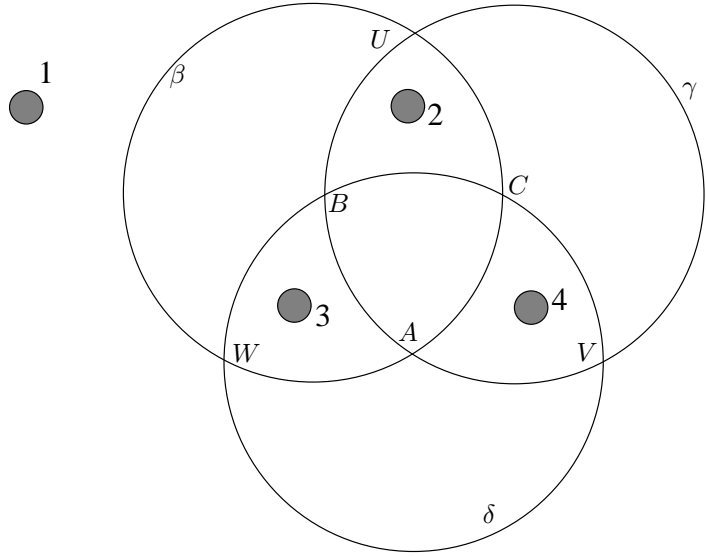}
\caption {A different view of Figure~\ref{fig:trio}. The four gray
  disks correspond to the four tubes from Figure~\ref{fig:trio}, and
  are marked accordingly.}
\label {fig:triad}
\end {center}
\end {figure}

We can construct similar Heegaard diagrams compatible with $L_0$ and $L_1$ as follows. The surface $\Sigma,$ the alpha curves and the basepoints remain the same. However, for $L_0$ we replace the beta curves by gamma curves $\gamma_i, \ i=1, \dots, n,$ while for $L_1$ we use delta curves $\delta_i, \ i=1, \dots, n.$ For $i < n,$ the curves $\gamma_i$ and $\delta_i$ are small isotopic translates of $\beta_i,$ such that they intersect $\beta_i$ in two points, and they also intersect each other in two points. For $i=n,$ we draw the curves $\gamma = \gamma_n$ and $\delta = \delta_n$ as in Figure~\ref{fig:trio}; see also Figure~\ref{fig:triad}, where the following intersection points are labelled:
$$ \beta \cap \gamma = \{ A, U \}, \ \ \gamma \cap \delta = \{B, V \}, \ \ \delta \cap \beta = \{C, W\}.$$

\begin {figure} \begin {center} \input {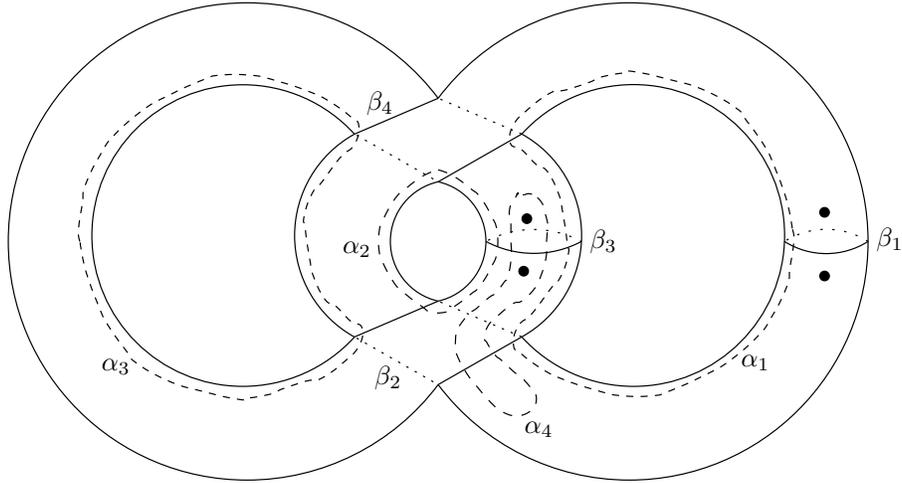}
\caption {A Heegaard diagram compatible with the Hopf link $L$, with $g=3, m=2$ 
and $n=4.$ The beta curves $\beta_2$ and $\beta = \beta_4$ are associated 
to the two crossings, $\beta_1$ to the distinguished edge, and $\beta_3$ 
is part of a ladybug.  There are three alpha curves associated to planar 
bounded regions and one, $\alpha_4$, which is part of a ladybug. One tip 
of $\alpha_4$ is stretched to achieve admissibility.} \label {fig:hopf} 
\end {center} \end {figure}

For the purpose of defining Floer homology, we need to ensure that the 
Heegaard diagrams for $L, L_0$ and $L_1$ constructed above are admissible 
in the sense of \cite[Definition 3.5]{OSlinks}. We achieve admissibility 
by stretching one tip of the alpha curve of each ladybug, and bringing it 
close to the basepoints associated to the distinguished edge $e.$ It is 
easy to see that the result is admissible; see Figure~\ref{fig:hopf} for 
an example. In that figure, to get the diagrams for $L_0$ and $L_1$, which 
are both the unknot, we replace $\beta = \beta_4$ by curves $\gamma$ and 
$\delta$, respectively, as in Figure~\ref{fig:trio}.

Now consider the tori
$$ \ta = \alpha_1 \times \dots \times \alpha_n, \ \  \tb = \beta_1 \times \dots \times \beta_n,$$
$$ \tc = \gamma_1 \times \dots \times \gamma_n, \ \
 \td = \delta_1 \times \dots \times \delta_n,$$
which we view as totally real 
submanifolds of the symmetric product $\sym^n(\wsigma).$ The Floer complex
$CF(\ta,
\tb)$ is the vector space freely generated by the intersection points
between $\ta$ and $\tb,$ and endowed with the differential
\begin{equation}
\label{eq:DefPartial}
 \del \xx = \sum_{\yy \in \ta \cap \tb} \sum_{\{\phi\in\hat\pi_2({\mathbf
    x},{\mathbf y})\big|
\mu(\phi)=1 \}}
\#\left(\frac{{\mathcal
      M}(\phi)}{\mathbb R}\right) {\mathbf y}.
\end{equation}
Here $\hat \pi_2({\mathbf x},{\mathbf y})$ denotes the space of
homology classes of Whitney disks connecting ${\mathbf x}$ to
${\mathbf y}$ in $\sym^n(\wsigma)$, ${\mathcal M}(\phi)$ denotes the
moduli space of pseudo-holomorphic representatives of $\phi$ (with
respect to a suitable almost complex structure as in~\cite{HolDisk}),
$\mu(\phi)$ denotes its formal dimension (Maslov index), and the $\#$
sign denotes the mod $2$ count of points in the (zero-dimensional)
moduli space. (We will henceforth use $\mu$ to denote Maslov index,
rather than the incidence number, as in Section~\ref{sec:khovanov}.)

The homology of $CF(\ta, \tb)$ is the Floer homology $HF(\ta, \tb).$ Up to a factor, this is the knot Floer homology of $L$:
$$ HF(\ta, \tb) \cong \hfk(L) \otimes V^{m-l},$$
where $V$ is a two-dimensional vector space as in Theorem~\ref{main}.

We can similarly take the Floer homology of $\ta$ and $\tc,$ or $\ta$ and $\td,$ and obtain
$$ HF(\ta, \tc) \cong \hfk(L_0) \otimes V^{m-l_0},$$
$$ HF(\ta, \td) \cong \hfk(L_1) \otimes V^{m-l_1}.$$

Therefore, the exact triangle from Theorem~\ref{main} can be written as
\begin {equation}
\label {ungraded}
\begin {CD}
HF(\ta, \td) @>(f_1)_*>> HF(\ta, \tb) @>(f_2)_*>>
HF(\ta,\tc) @>(f_3)_*>> HF(\ta, \td) 
\end {CD}
\end {equation}

The maps $(f_i)_* \ (i=1,2,3)$ from the triangle \eqref{ungraded} are all 
induced by chain maps $f_i$ between the corresponding Floer complexes. To 
define the maps $f_i$, let us first recall the definition of the usual 
triangle maps appearing in Floer theory. Given totally real submanifolds 
$T_1, T_2, T_3$ in a symplectic manifold (satisfying several technical 
conditions which will hold in our situations), 
there is a chain map 
$$ CF(T_1, T_2)\otimes CF(T_2, T_3) \to CF(T_1, T_3),$$ 
defined by counting pseudo-holomorphic triangles.
In particular, given an intersection point ${\mathbf z}\in T_2\cap T_3$
which is a cycle when viewed as an element of $CF(T_2,T_3)$, we have a chain map
$$ F_{\mathbf z} (\xx) 
= \sum_{\yy \in T_1 \cap T_3} \sum_{\{\phi\in \hat\pi_2({\mathbf
    x},{\mathbf z}, {\mathbf y})\big|
\mu(\phi)=0 \}}
\#\left({\mathcal
      M}(\phi)\right) {\mathbf y}. $$
Here $\hat \pi_2({\mathbf x},{\mathbf z}, {\mathbf y})$ denotes the space of 
homology classes of triangles with edges on $T_1, T_2, T_3$ and vertices 
$\xx, {\mathbf z}$ and $\yy,$ respectively (in clockwise order), $\mu$ is 
the Maslov index, and $\#\left({\mathcal M}(\phi)\right)$ the number of 
their pseudo-holomorphic representatives.

Going back to our set-up, whenever we have two isotopic curves $\eta$ 
and $\eta'$ on the surface $\Sigma$ such that they intersect in exactly 
two points, we will denote by $M_{\eta \eta'}\in \eta \cap \eta'$ the top 
degree generator of $CF(\eta, \eta').$ Given one of the intersection points in 
Figure~\ref{fig:triad}, for example $A 
\in \beta \cap \gamma,$ we
obtain a corresponding intersection point 
in $\tb \cap \tc$ by adjoining to $A$ the top degree intersection points
$M_{\beta_i\gamma_i} \in \beta_i \cap \gamma_i.$ We denote the resulting
generators by the respective lowercase letters in bold:
$$ \mathbf{a} = M_{\beta_1\gamma_1} \times M_{\beta_2\gamma_2}\times \dots
\times M_{\beta_{n-1}\gamma_{n-1}}\times A \in
\tb \cap  \tc;$$
$$ \mathbf{b} = M_{\gamma_1\delta_1}\times M_{\gamma_2\delta_2} \times
\dots \times M_{\gamma_{n-1}\delta_{n-1}} \times B
\in \tc \cap \td; $$
$$ \mathbf{c} = M_{\delta_1\beta_1}\times M_{\delta_2\beta_2} \times
\dots \times M_{\delta_{n-1}\beta_{n-1}} \times C
\in \td \cap \tb; $$ 
$$ \mathbf{u} = M_{\beta_1\gamma_1} \times M_{\beta_2\gamma_2}\times \dots
\times M_{\beta_{n-1}\gamma_{n-1}}\times U \in
\tb \cap  \tc;$$ 
$$ \mathbf{v} = M_{\gamma_1\delta_1}\times M_{\gamma_2\delta_2} \times
\dots \times M_{\gamma_{n-1}\delta_{n-1}} \times V
\in \tc \cap \td; $$              
$$ \mathbf{w} = M_{\delta_1\beta_1}\times M_{\delta_2\beta_2} \times  
\dots \times M_{\delta_{n-1}\beta_{n-1}} \times W
\in \td \cap \tb. $$

The chain maps $f_i$ giving rise to \eqref{ungraded} are then defined to 
be the sums 
$$ f_1 = F_{\mathbf{c}}+F_{\mathbf{w}} : CF(\ta, \td) \to CF(\ta, \tb);$$
$$ f_2 = F_{\mathbf{a}}+F_{\mathbf{u}} : CF(\ta, \tb) \to CF(\ta, \tc);$$
$$ f_3 = F_{\mathbf{b}}+F_{\mathbf{v}} : CF(\ta, \tc) \to CF(\ta, \td).$$

\subsection {Periodic domains}
\label {sec:per}

Let us apply the discussion in Section~\ref{sec:periodic} to the setting of Section~\ref{sec:heegaard}. 

Note that $(\Sigma; \aalpha, \bbeta)$, for example, is a Heegaard diagram for $S^3$, hence the alpha and the beta curves span all of $H_1(\Sigma; \zz).$ Applying Corollary~\ref{pery} we deduce that 
\begin {equation}
\label {carnation}
\per_{\alpha, \beta} = \per_{\alpha} + \per_{\beta}.
\end {equation}
Similarly, we have $\per_{\alpha, \gamma} = \per_{\alpha} + \per_{\gamma}$ and $\per_{\alpha, \delta} = \per_{\alpha} + \per_{\delta}.$

The situation for $\per_{\beta, \gamma}$ is different. Before analyzing it, let us first understand the components of $\Sigma - (\cup \beta_i)$, which span $\per_{\beta}$, in detail. Their number is $m$, which equals either $l$ or $l+1,$ according to whether the two strands of $L$ meeting at $c$ are on different link components, or on the same link component. Let $K_1, \dots, K_l$ be the connected components of $L$, such that $K_l$ is the one containing the edge $e.$ If $m=l,$ then each $K_i$ corresponds to a unique component $\dom_i^\beta$ of $\Sigma - (\cup \beta_i)$, which lies in a neighborhood of $K_i$ (when $\Sigma$ is viewed as the boundary of a neighborhood of $L$). If $m = l+1,$ then for $i < l$, each $K_i$ corresponds again to some $\dom^{\beta}_i,$ but in the neighborhood of $K_l$ there are now two components of $\Sigma - (\cup \beta_i)$, which we denote by $\dom^{\beta}_l$ and $\dom^{\beta}_{l+1}$, such that $\dom^{\beta}_l$ is the one whose boundary contains the curve $\beta = \beta_n.$

\begin {figure}
\begin {center}
\begin{picture}(0,0)%
\includegraphics{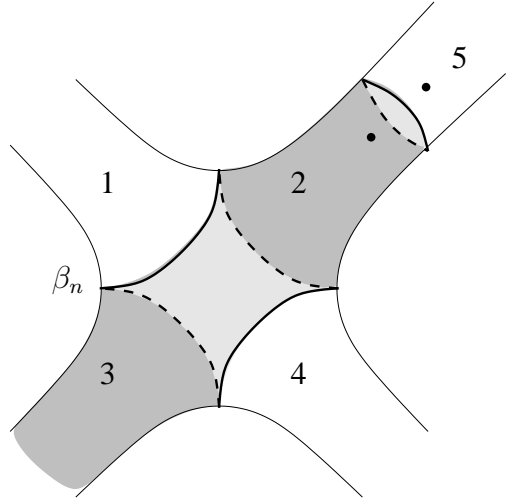}%
\end{picture}%
\setlength{\unitlength}{1973sp}%
\begingroup\makeatletter\ifx\SetFigFont\undefined%
\gdef\SetFigFont#1#2#3#4#5{%
  \reset@font\fontsize{#1}{#2pt}%
  \fontfamily{#3}\fontseries{#4}\fontshape{#5}%
  \selectfont}%
\fi\endgroup%
\begin{picture}(6249,6249)(1264,-4798)
\put(1801,-2161){\makebox(0,0)[lb]{\smash{{\SetFigFont{12}{14.4}{\rmdefault}{\mddefault}{\updefault}{\color[rgb]{0,0,0}$\beta_n$}%
}}}}
\end{picture}%

\caption{We illustrate here Equation~\eqref{sbeta}.
The component $\dom^{\beta}_l$, which gives a homological relation
between $\beta_n$ and other $\beta$-curves, is shaded. There are two
cases: when $m=\ell$, the region labelled here by $5$ is included in
$\dom^{\beta}_l$. Otherwise, when $m=\ell+1$, $\dom^{\beta}_l$
terminates in a different meridinal $\beta$-circle. In either case,
the boundary of $\dom^{\beta}_l$ consists of $\beta$-circles, and it
contains $\beta_n$ with multiplicity one.}
\label{fig:PeriodicDomain}
\end {center}
\end {figure}

Note that, regardless of whether $m=l$ or $m=l+1$, the component
$\dom^{\beta}_l$ contains the curve $\beta_n$ with multiplicity $\pm
1$ (see Figure~\ref{fig:PeriodicDomain}). This means that the class
$[\beta_n] \in S_{\beta} \subset H_1(\Sigma; \zz)$ is in the span of
the other beta curves. In other words,
\begin {equation}
\label {sbeta}
S_{\beta} = \text{ Span } (\beta_1, \dots, \beta_{n-1}).
\end {equation}

Similar remarks apply to $\Sigma - (\cup \gamma_i)$ and $\Sigma - (\cup \delta_i)$. Their components are denoted $\dom^{\gamma}_i$ and $\dom^{\delta}_i$, respectively, for $i=1, \dots, m.$ Recall that for each $i=1, \dots, n-1$, the curves $\beta_i$,  $\gamma_i$ and $\delta_i$ are isotopic. Therefore, Equation~\eqref{sbeta}, together with its analogs for the beta and gamma curves, implies that
\begin {equation}
\label {sbcd}
S_{\beta} = S_{\gamma} = S_{\delta}.
\end {equation}

For each $j=1, \dots, n-1$, the curves $\beta_j$ and $\gamma_j$ are separated by two thin bigons in $\Sigma.$ The difference of these bigons is a periodic domain $\dom^{\beta, \gamma}_j,$ with boundary $\beta_j - \gamma_j.$ Equation~\eqref{sbcd} implies that $\rk (S_{\beta, \gamma}) = \rk(S_{\beta}) = g,$ so from Lemma~\ref{perab} we deduce that $\rk (\per_{\beta, \gamma}) = 2n+1 -g = n+m.$ In fact, it is not hard to check that the following is true:
\begin {lemma}
\label {freud}
The domains $\dom^{\beta}_i, \dom^{\gamma}_i \ (i=1, \dots, m)$ and $\dom^{\beta, \gamma}_j \ (j=1, \dots, n-1)$ span the group $\per_{\beta, \gamma}.$
\end {lemma}
 
Note that we gave a set of $2m + n-1$ generators for the group $\per_{\beta, \gamma}$ of rank $n+m.$ There are indeed $m-1$ independent relations between these generators, namely for each of the $m-1$ components $K_i$ of $L$ (or $L_0$) not containing either of the strands intersecting at $c$, the difference $\dom^{\beta}_i - \dom^{\gamma}_i$ can also be written as a sum of some domains $\dom^{\beta, \gamma}_j$ (corresponding to the crossings on $K_i$).

Next, let us look at the triply periodic domains with boundary on the alpha, beta, and gamma curves.

\begin {lemma}
\label {jung}
We have $\per_{\alpha, \beta, \gamma}  = \per_{\alpha} + \per_{\beta, \gamma}.$
\end {lemma}

\noindent\textbf{Proof.} Let $\dom$ be any triply periodic domain in $\per_{\alpha, \beta, \gamma}.$ If the curve $\gamma_n$ appears (with nonzero multiplicity) in the boundary of $\dom,$ by the analog of \eqref{sbeta} for gamma curves we can subtract some domain in $\per_{\gamma} \subset \per_{\beta, \gamma}$ from $\dom$ and obtain a new domain, in which the multiple of $\gamma_n$ from $\del \dom$ was traded for a combination of the other gamma curves $\gamma_1, \dots, \gamma_{n-1}.$ Next, whenever we have some curve $\gamma_j$ in the boundary ($j < n$), we can add the corresponding domain $ \dom^{\beta, \gamma}_j \in \per_{\beta} \subset \per_{\beta, \gamma}$ to trade it for a beta curve. Thus we arrive at a domain in $\per_{\alpha, \beta}$ and the conclusion follows from Equation~\eqref{carnation}. $\hfill \fin$

\medskip

Note that Lemma~\ref{freud} has straightforward analogs about the structure of the groups $\per_{\gamma, \delta}$ and $ \per_{\delta, \beta}.$ Similarly, Lemma~\ref{jung} has straightforward analogs about the structure of the groups $\per_{\alpha, \gamma, \delta}$ and $ \per_{\alpha, \delta, \beta}.$

\subsection {The relative $\delta$-grading}

Pick $\xx, \yy \in \ta \cap \tb$. Let $\pi_2(\xx,\yy)$ be the space of  homology classes of Whitney disks in $\sym^n(\Sigma)$ connecting $\xx$ and $\yy.$ (Recall that $\hat \pi_2(\xx,\yy)$ is the corresponding space in $\sym^n(\wsigma).$) Since $(\Sigma, \alpha_1, \dots,\alpha_n, \beta_1, \dots,\beta_n)$ is a Heegaard diagram for $S^3$, we have $\pi_2(\xx,\yy) \neq \emptyset$ for any $\xx$ and $\yy.$ Note that $\pi_2(\xx, \xx)$ can be identified with the group of periodic domains $\per_{\alpha, \beta}.$

Every class $\phi \in \pi_2(\xx,\yy)$ has a Maslov index $\mu(\phi)
\in \zz.$ In the usual construction of knot Floer homology, the extra
basepoints on the Heegaard surface $\Sigma$ are of two types: half of
them are denoted $w_j$ and the other half $z_j,$ with $j=1, \cdots,
m+1$, such that every connected component of $\Sigma - \cup \alpha_i $
or $\Sigma - \cup \beta_i$ contains exactly one of the $w_j$ and one of the $z_k$. Let
$W(\phi)$ and $Z(\phi)$ be the intersection numbers of $\phi$ with the
union of all $\{ w_j \} \times \sym^{n-1}(\Sigma)$ and the union of
all $\{ z_j \} \times \sym^{n-1}(\Sigma)$, respectively. Thus $ \hat
\pi_2(\xx,\yy)$ is the space of classes $\phi$ with $W(\phi) =
Z(\phi)=0.$
      
The difference in the Maslov grading $H$ (denoted $i$ in the
introduction) between $\xx$ and $\yy$ can be calculated by picking
some $\phi \in \pi_2(\xx,\yy)$ and applying the formula
$$ H(\xx) -H(\yy) = \mu(\phi) - 2W(\phi).$$

Similarly, the difference in the Alexander grading $A$ (denoted $j$ in the introduction) is
$$ A(\xx) -A(\yy) = Z(\phi) - W(\phi).$$

Setting $P(\phi) = Z(\phi) + W(\phi)$, the difference in the grading $\delta = A - H$ is then
$$ \delta(\xx) -\delta(\yy) = P(\phi) - \mu(\phi).$$
 
Therefore, if we limit ourselves to considering the $\delta$ grading,
there is no difference between the two types of basepoints. This
explains why we have not distinguished between them in
Section~\ref{sec:heegaard}, and we will not distinguish between them
from now on either.

Observe that the relative $\delta$ grading is well-defined, i.e. we
have $\mu(\phi)-P(\phi)= \mu(\phi') - P(\phi')$ for any $\phi, \phi' \in \pi_2(\xx, \yy)$. Indeed, because $\mu$ and $P$ are additive under concatenation, it suffices to prove that $\mu(\phi)-P(\phi)=0$ for any $\phi \in \pi_2(\xx,
\xx)= \per_{\alpha \beta}$. By Equation \eqref{carnation}, the group $\per_{\alpha \beta}$ is generated by the connected components of $\Sigma - \cup
\alpha_i $ and $\Sigma - \cup \beta_i$. Each such component has
$\mu(\phi)=P(\phi)=2,$ so the relative $\delta$ grading is well-defined.

\begin {lemma}
\label {reldelta}
The chain maps $f_1, f_2, f_3$ that induce the triangle  
\eqref{ungraded} preserve the relative $\delta$ grading.
\end{lemma}

\noindent\textbf{Proof.} First, observe that a triangle map such as $ 
F_{\mathbf{a}} : CF(\ta, \tb) \to CF(\ta, \tc)$ preserves the relative 
$\delta$ grading. In other words, we need to show that adding a triply periodic domain $\dom \in \per_{\alpha, \beta, \gamma}$ to a class $\phi \in \hat \pi_2(\xx, \mathbf{a}, \yy)$ does not change the quantity $\mu(\phi)-P(\phi).$ By Lemmas~\ref{freud} and \ref{jung}, it suffices to show that the classes of the domains $\dom^{\alpha}_i, \dom^{\beta}_i, \dom^{\gamma}_i$ and $\dom^{\beta, \gamma}_j$ all have $\mu = P.$ Indeed, for  $\dom^{\alpha}_i, \dom^{\beta}_i$ and $\dom^{\gamma}_i$ this is the argument in the paragraph before Lemma~\ref{reldelta}, while for each  $\dom^{\beta, \gamma}_j \ (j=1, \dots, n-1)$ we have $\phi = P =0.$

\begin {figure}
\begin {center}
\input {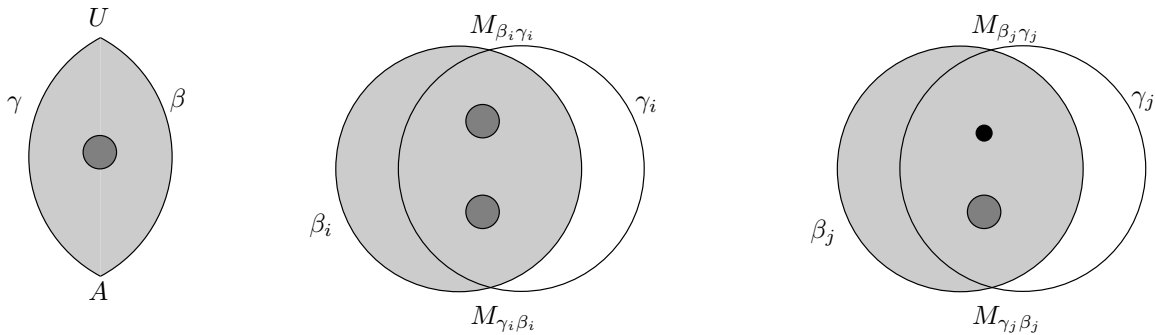}
\caption {A relative homology class $\phi  \in \pi_2(\aa,\uu)$ of Maslov index $1,$ with one basepoint. The grey disks here denote tubes, whereas the small
black dot in the rightmost picture denotes a basepoint.}
\label {fig:au}
\end {center}
\end {figure}

Next, in order to show that $ f_2 = F_{\mathbf{a}}+F_{\mathbf{u}}$
preserves the relative $\delta$-grading, we exhibit a class $\phi \in
\pi_2(\aa,\uu)$ with $\mu(\phi)=P(\phi).$ In Figure~\ref{fig:triad}
there is a bigon relating $A$ and $U$ which is connected by the tube
numbered 2 to the rest of the Heegaard diagram. This bigon is also
shown on the left in Figure~\ref{fig:au}. Following the tube, we
encounter several disks (or possibly none) bounded by beta circles as
in the middle of Figure~\ref{fig:au}, until we find a disk as on the
right of Figure~\ref{fig:au}. Lipshitz's formula for the Maslov index
\cite{Li} says that $\mu(\phi)$ can be computed as the sum of the
Euler measure $e(\phi)$ and a vertex multiplicity $n(\phi).$ (We refer
to \cite{Li} for the definitions.) The punctured bigon on the left of
Figure~\ref{fig:au} contributes $-\frac{1}{2}$ to $e(\phi)$ and
$\frac{1}{2}$ to $n(\phi)$, each middle disk $-1$ to $e(\phi)$ and $1$
to $n(\phi)$, and the disk on the right $0$ to $e(\phi)$ and $1$ to
$n(\phi)$. Thus $\mu(\phi)=P(\phi)=1.$

The arguments for $f_1$ and $f_3$ are similar. $\hfill \fin$

\subsection {The absolute $\delta$-grading}
\label {sec:absdelta}
The generators $\xx\in \ta \cap \tb$ are of two kinds. They all consist of 
$n$-tuples of points in $\Sigma$, one on each alpha curve and on each 
beta curve. If for each ladybug (consisting of a pair of curves $\alpha_i$ 
and $\beta_i$), $\xx$ contains one of the two points in $ \alpha_i \cap 
\beta_i$, we call the generator $\xx$ \textit{Kauffman}. Otherwise, it is 
called \textit{non-Kauffman}. Note that, if we hadn't had to stretch the 
alpha curves on the ladybugs to achieve admissibility, all generators 
would have been Kauffman.

Every $\xx\in \ta \cap \tb$ has an absolute $\delta$-grading $\delta(\xx) 
\in \frac{1}{2} \zz.$ We will explain now a simple formula for 
$\delta(\xx)$ when $\xx$ is Kauffman.

Consider the regions $A_0, A_1, A_2, \dots, A_{k+1}$ as in the second
paragraph after the statement of Theorem~\ref{main}. Each of the $k$
crossings in $D$ is on the boundary of four regions. A \textit{state},
cf. \cite{Ka}, is an assignment which associates to each crossing one
of the four incoming quadrants, such that the quadrants associated to
distinct vertices are in distinct regions, and none are corners of the
regions $A_0$ or $A_1.$

One can associate a monomial to each state such that as we sum all
these monomials we obtain the Alexander polynomial of the link $L$,
\cite{Ka}.  Therefore, if the Alexander polynomial $\Delta_L(q)$ is
nonzero (or, in particular, if $\Delta_L(-1) = \det(L) \neq 0$), then
there must be at least one state.

To each Kauffman generator $\xx$ we can associate a state in an
natural way: at each crossing $c$ the corresponding beta curve
intersects exactly one of the alpha curves of the neighboring regions
in a point of $\xx$, and the quadrant in that region is the one we
associate to $c.$ In \cite{OSalt}, the Maslov and Alexander gradings
of Kauffman generators are calculated in terms of their states;
compare also~\cite{OSS}.

\begin {figure}
\begin {center}
\input {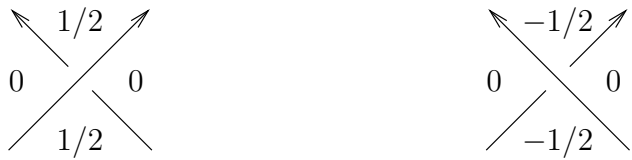}
\caption {Contributions to the $\delta$-grading.}
\label {fig:contrib}
\end {center}
\end {figure}

For our purposes, it suffices to know how to compute the $\delta$-grading. If $\xx$ is Kauffman and $c$ is a crossing in $D$, we let $\delta(\xx,c) \in \{0, \pm 1/2 \}$ be the quantity from Figure~\ref{fig:contrib}, chosen according to which quadrant at $c$ appears in the state of $\xx.$ Then:
\begin {equation}
\label {delta}
\delta(\xx) = \sum_c  \delta(\xx,c).
\end {equation}

A similar discussion applies to the diagrams $D_0$ and $D_1$ of the resolutions $L_0$ and $L_1$, respectively, except that in those cases there is no contribution from the resolved crossing $c_0.$ 

Note that the $\delta$-grading of a Kauffman generator $\xx$ does not depend on which of the two intersection points between the two curves of a ladybug appears in $\xx.$ 

\begin {lemma}
\label {eins}
Suppose that $c_0$ is a positive crossing in $D$ (so that $L_0$ is the oriented resolution $L_v$) and that $\det(L_0) \neq 0.$ Then the map $f_2 : CF(\ta, \tb) \to CF(\ta, \tc)$ decreases $\delta$-grading by $1/2.$
\end {lemma}

\noindent\textbf{Proof.} By Lemma~\ref{reldelta}, we already know that 
$f_2$ preserves the relative $\delta$-grading. Thus, it suffices to 
exhibit two generators $\xx\in \ta\cap\tb$ and $\yy\in \ta\cap\tc$ with 
$\delta(\xx) - \delta(\yy) =1/2,$ and such that there exists a 
holomorphic triangle of index zero in $\hat \pi_2({\mathbf x},{\mathbf a}, 
{\mathbf y}).$

Since $\det(L_0) \neq 0,$ the diagram $D_0$ has at least one Kauffman 
generator $\yy.$ There is a corresponding Kauffman generator $\xx\in 
\ta\cap\tb$, such that each $y_i \in \gamma_i \cap \yy, \ (i < n)$ is 
close to some $x_i \in \beta_i \cap \xx$ (they are related by the isotopy 
between $ \gamma_i$ and $\beta_i$), while $x_n \in \beta \cap \xx$ and 
$y_n \in \gamma \cap \yy$ are two vertices of the shaded triangle in 
Figure~\ref{fig:trio} with the third vertex at $A.$ That shaded triangle, 
coupled with the small triangles with vertices at $\xx_i, \yy_i,$ and 
$M_{\beta_i, \gamma_i} $ for $i =1, \dots, n-1$, gives the desired 
holomorphic triangle in $\sym^n(\wsigma).$ To check that $\delta(\xx) - 
\delta(\yy) =1/2,$ note that in formula \eqref{delta} the contributions 
to $\delta(\xx)$ and $\delta(\yy)$ from each crossing are the same, except 
that there is an extra contribution of $1/2$ to $\delta(\xx)$ coming from 
$c_0. \ \hfill \fin$

\begin {lemma}
\label{zwei}
Suppose that $c_0$ is a positive crossing in $D$ (so that $L_1$ is the
unoriented resolution $L_h$) and that $\det(L_1) \neq 0.$ Then the map
$f_1 : CF(\ta, \td) \to CF(\ta, \tb)$ shifts $\delta$-grading by
$e/2$, where $e$ is as in the statement of Lemma~\ref{signatures}.
\end {lemma}

\noindent\textbf{Proof.} By Lemma~\ref{reldelta}, we already know that 
$f_1$ preserves the relative $\delta$-grading. Again, it suffices to 
exhibit two generators $\xx\in \ta\cap\tb$ and $\yy\in \ta\cap\td$ with 
$\delta(\xx) - \delta(\yy) =e/2,$ and such that there exists a holomorphic 
triangle of index zero in $\hat \pi_2({\mathbf y},{\mathbf w}, {\mathbf 
x}).$

Since $\det(L_1) \neq 0,$ we can pick a Kauffman generator  $\yy\in 
\ta\cap\td$. As in the proof of Lemma~\ref{eins}, there is a corresponding 
Kauffman generator $\xx\in \ta\cap\tb$ and a holomorphic triangle of index 
zero as desired, consisting of $n-1$ small triangles with one vertex at 
$M_{\delta_i, \beta_i} $ for $i =1, \dots, n-1$, and the shaded triangle in Figure~\ref{fig:trio} with one vertex at $W.$  

To check that $\delta(\xx) - \delta(\yy) =e/2,$ let $n_+$ be the number of positive crossings in $D$ (excluding $c_0$) which change sign in $D_1.$ At each such crossing $c,$ we have:
$$ \delta(\xx, c) = \delta(\yy, c) + 1/2.$$

Let also $n_-$ be the number of negative crossings in $D$ which change sign in $D_1.$ At each such crossing $c,$ we have:
$$ \delta(\xx, c) = \delta(\yy, c) - 1/2.$$

Therefore,
$$ \delta(\xx) - \delta(\yy) = (n_+ - n_-)/2 = e/2.$$ $\hfill \fin$
 
\begin {proposition}
\label {absgrad}
Let $L$ be a link and $L_0, L_1$ its two resolutions at a crossing as in Figure~\ref{fig:skein}. Assume that $\det(L_0), \det(L_1) > 0$ and $\det(L) = \det(L_0) + \det(L_1).$ Then, two of the three maps in the exact triangle from Theorem~\ref{main} behave as follows with respect to the $\delta$-grading:
$$
 \hfk_{*-\frac{\sigma(L_1)}{2}}(L_1) \otimes V^{m-l_1} \To \hfk_{*-\frac{\sigma(L)}{2}}(L) 
 \otimes V^{m-l_0} \To \hfk_{*-\frac{\sigma(L_0)}{2}}(L_0) \otimes V^{m-l} ,$$
where $V$ denotes a two-dimensional vector space over $\ff$, in grading zero.
\end {proposition}

\noindent\textbf{Proof.} When the given crossing in $L$ is positive, this follows from \eqref{ungraded}, together with the results of Lemmas~\ref{signatures}, \ref{eins}, and \ref{zwei}. The case when the crossing is negative is similar. $\hfill \fin$
\medskip

\noindent\textbf{Proof of Theorem~\ref{floer}.} Using Proposition~\ref{absgrad}, we can argue in the same way as in the proof of Theorem~\ref{khovanov}. Note that we do not have to know the change in the absolute $\delta$-grading under the third map $(f_3)_*: HF(\ta, \tc) \to HF(\ta, \td)$ in the exact triangle. Indeed, recall that the Euler characteristic of $\hfk$ is (up to a factor) the Alexander polynomial, which evaluated at $-1$ gives the determinant of the link. If we know that $L_0$ and $L_1$ are Floer homologically $\sigma$-thin and we want to show the same for $L$, the fact that $\det(L) = \det(L_0) + \det(L_1)$ together with the ungraded triangle implies that 
$$  \rk \ (\hfk(L) \otimes V^{m-l}) = \rk \ (\hfk(L_0) \otimes V^{m-l_0}) + \rk \ (\hfk(L_1) \otimes V^{m-l_1}).$$
Hence $(f_3)_*=0$, and the inductive step goes through. $\hfill \fin$

\end{document}